\documentclass[11pt]{article}

\usepackage{amssymb}

\newcommand{\newc}{\newcommand}

\newc{\eqnoset}{\setcounter{equation}{0}}

\newcommand{\mref}[1]{(\ref{#1})}
\newcommand{\reflemm}[1]{Lemma~\ref{#1}}
\newcommand{\refrem}[1]{Remark~\ref{#1}}
\newcommand{\reftheo}[1]{Theorem~\ref{#1}}

\newcommand{\refcoro}[1]{Corollary~\ref{#1}}
\newcommand{\refprop}[1]{Proposition~\ref{#1}}
\newcommand{\refsec}[1]{Section~\ref{#1}}

\newcommand{\beq}{\begin{equation}}
\newcommand{\eeq}{\end{equation}}
\newcommand{\beqno}[1]{\begin{equation}\label{#1}}

\newcommand{\barr}{\begin{array}}
\newcommand{\earr}{\end{array}}

\newc{\bearr}{\begin{eqnarray*}}
\newc{\eearr}{\end{eqnarray*}}

\newc{\bearrno}[1]{\begin{eqnarray}\label{#1}}
\newc{\eearrno}{\end{eqnarray}}

\newc{\non}{\nonumber}
\newc{\nol}{\nonumber\nl}

\newcommand{\bdes}{\begin{description}}
\newcommand{\edes}{\end{description}}
\newc{\benu}{\begin{enumerate}}
\newc{\eenu}{\end{enumerate}}
\newc{\btab}{\begin{tabular}}
\newc{\etab}{\end{tabular}}

\newtheorem{theorem}{Theorem}[section]
\newtheorem{defi}[theorem]{Definition}
\newtheorem{lemma}[theorem]{Lemma}
\newtheorem{rem}[theorem]{Remark}
\newtheorem{exam}[theorem]{Example}
\newtheorem{propo}[theorem]{Proposition}
\newtheorem{corol}[theorem]{Corollary}

\newcommand{\btheo}[1]{\begin{theorem}\label{#1}}
\newc{\brem}[1]{\begin{rem}\label{#1}\em}
\newc{\bexam}[1]{\begin{exam}\label{#1}\em}
\newc{\bdefi}[1]{\begin{defi}\label{#1}}
\newcommand{\blemm}[1]{\begin{lemma}\label{#1}}
\newcommand{\bprop}[1]{\begin{propo}\label{#1}}
\newcommand{\bcoro}[1]{\begin{corol}\label{#1}}
\newcommand{\etheo}{\end{theorem}}
\newcommand{\elemm}{\end{lemma}}
\newcommand{\eprop}{\end{propo}}
\newcommand{\ecoro}{\end{corol}}
\newc{\erem}{\end{rem}}
\newc{\eexam}{\end{exam}}
\newc{\edefi}{\end{defi}}

\newc{\rmk}[1]{{\bf REMARK #1: }}
\newc{\DN}[1]{{\bf DEFINITION #1: }}

\newcommand{\bproof}{{\bf Proof:~~}}
\newc{\eproof}{{\vrule height8pt width5pt depth0pt}\vspace{3mm}}

\newc{\bfrac}[2]{\dspl{\frac{#1}{#2}}}

\newc{\nid}{\noindent}

\newcommand{\dspl}{\displaystyle}
\newc{\grad}{\nabla}
\newc{\Div}{\mbox{div}}
\newc{\pdt}[1]{\dspl{\frac{\partial{#1}}{\partial t}}}
\newc{\pdn}[1]{\dspl{\frac{\partial{#1}}{\partial \nu}}}
\newc{\pdNi}[1]{\dspl{\frac{\partial{#1}}{\partial \mathcal{N}_i}}}
\newc{\pD}[2]{\dspl{\frac{\partial{#1}}{\partial #2}}}
\newc{\dt}{\dspl{\frac{d}{dt}}}
\newc{\bdry}[1]{\mbox{$\partial #1$}}
\newc{\sgn}{\mbox{sign}}

\newc{\Hess}[1]{\frac{\partial^2 #1}{\pdh z_i \pdh z_j}}
\newc{\hess}[1]{\partial^2 #1/\pdh z_i \pdh z_j}

\newc{\ag}{\alpha}
\newc{\bg}{\beta}
\newc{\cg}{\gamma}\newc{\Cg}{\Gamma}
\newc{\dg}{\delta}\newc{\Dg}{\Delta}
\newc{\eg}{\varepsilon}
\newc{\zg}{\zeta}
\newc{\thg}{\theta}
\newc{\llg}{\lambda}\newc{\LLg}{\Lambda}
\newc{\kg}{\kappa}
\newc{\rg}{\rho}
\newc{\sg}{\sigma}\newc{\Sg}{\Sigma}
\newc{\tg}{\tau}
\newc{\fg}{\phi}\newc{\Fg}{\Phi}
\newc{\vfg}{\varphi}
\newc{\og}{\omega}\newc{\Og}{\Omega}
\newc{\pdh}{\partial}

\newc{\ccG}{{\cal G}}

\newc{\ii}[1]{\int_{#1}}
\newc{\iidx}[2]{{\dspl\int_{#1}~#2~dx}}
\newc{\bii}[1]{{\dspl \ii{#1} }}
\newc{\biii}[2]{{\dspl \iii{#1}{#2} }}
\newc{\su}[2]{\sum_{#1}^{#2}}
\newc{\bsu}[2]{{\dspl \su{#1}{#2} }}

\newc{\biiom}[1]{{\dspl\int_{\bdrom}~ #1 ~d\sg}}
\newc{\io}[1]{{\dspl\int_{\Og}~ #1 ~dx}}
\newc{\bio}[1]{{\dspl\int_{\bdrom}~ #1 ~d\sg}}
\newc{\bsir}{\bsu{i=1}{r}}
\newc{\bsim}{\bsu{i=1}{m}}

\newc{\iibr}[2]{\iidx{\bprw{#1}}{#2}}
\newc{\Intbr}[1]{\iibr{R}{#1}}
\newc{\intbr}[1]{\iibr{\rg}{#1}}
\newc{\intt}[3]{\int_{#1}^{#2}\int_\Og~#3~dxdt}
\newc{\itQ}[2]{\dspl{\int\hspace{-2.5mm}\int_{#1}~#2~dz}}
\newc{\mitQ}[2]{\dspl{\rule[1mm]{4mm}{.3mm}\hspace{-5.3mm}\int\hspace{-2.5mm}\int_{#1}~#2~dz}}
\newc{\mitQQ}[3]{\dspl{\rule[1mm]{4mm}{.3mm}\hspace{-5.3mm}\int\hspace{-2.5mm}\int_{#1}~#2~#3}}

\newc{\mitx}[2]{\dspl{\rule[1mm]{3mm}{.3mm}\hspace{-4mm}\int_{#1}~#2~dx}}
\newc{\mitmu}[2]{\dspl{\rule[1mm]{3mm}{.3mm}\hspace{-4mm}\int_{#1}~#2~d\mu}}
\newc{\iidmu}[2]{{\dspl\int_{#1}~#2~d\mu}}

\newc{\iidm}[3]{{\dspl\int_{#1}~#2~d #3}}

\newc{\itQmu}[2]{\dspl{\int\hspace{-2.5mm}\int_{#1}~#2~d\mu}}
\newc{\mitQmu}[2]{\dspl{\rule[1mm]{4mm}{.3mm}\hspace{-5.3mm}\int\hspace{-2.5mm}\int_{#1}~#2~d\mu}}

\newc{\mitQq}[2]{\dspl{\rule[1mm]{4mm}{.3mm}\hspace{-5.3mm}\int\hspace{-2.5mm}\int_{#1}~#2~d\bar{z}}}
\newc{\itQq}[2]{\dspl{\int\hspace{-2.5mm}\int_{#1}~#2~d\bar{z}}}

\newc{\pder}[2]{\dspl{\frac{\partial #1}{\partial #2}}}

\newc{\bdrom}{\bdry{\Og}}

\newc{\bilhom}{\mbox{Bil}(\mbox{Hom}(\RR^{nm},\RR^{nm}))}
\newc{\VV}[1]{{V(Q_{#1})}}

\newc{\ccA}{{\mathcal A}}
\newc{\ccB}{{\mathcal B}}
\newc{\ccC}{{\mathcal C}}
\newc{\ccD}{{\mathcal D}}
\newc{\ccE}{{\mathcal E}}
\newc{\ccH}{\mathcal{H}}
\newc{\ccF}{\mathcal{F}}
\newc{\ccI}{{\mathcal I}}
\newc{\ccJ}{{\mathcal J}}
\newc{\ccK}{{\mathcal K}}
\newc{\ccP}{{\mathcal P}}
\newc{\ccQ}{{\mathcal Q}}
\newc{\ccR}{{\mathcal R}}
\newc{\ccS}{{\mathcal S}}
\newc{\ccT}{{\mathcal T}}
\newc{\ccX}{{\mathcal X}}
\newc{\ccY}{{\mathcal Y}}
\newc{\ccZ}{{\mathcal Z}}

\newc{\bb}[1]{{\mathbf #1}}

\newc{\myprod}[1]{\langle #1 \rangle}
\newc{\mypar}[1]{\left( #1 \right)}
\newc{\BLLg}{\mathbf{\LLg}}

\newc{\mA}{\mathbf{A}}
\newc{\mB}{\mathbf{B}}
\newc{\mC}{\mathbf{C}}
\newc{\mD}{\mathbf{D}}
\newc{\mE}{\mathbf{E}}
\newc{\mF}{\mathbf{F}}
\newc{\mJ}{\mathbf{J}}
\newc{\mG}{\mathbf{G}}
\newc{\mP}{\mathbf{P}}
\newc{\mR}{\mathbf{R}}
\newc{\mQ}{\mathbf{Q}}
\newc{\mX}{\mathbf{X}}
\newc{\muu}{\mathbf{u}}
\newc{\mvv}{\mathbf{v}}

\newc{\mllg}{\mathbb{\lambda}}
\newc{\mLLg}{\mathbf{\LLg}}

\newc{\lspn}[2]{\mbox{$\| #1\|_{\Lsp{#2}}$}}
\newc{\Lpn}[2]{\mbox{$\| #1\|_{#2}$}}
\newc{\Hn}[1]{\mbox{$\| #1\|_{H^1(\Og)}$}}

\newc{\mynorm}[2]{\| #1\|_{#2}}

\newcommand{\RR}{{\rm I\kern -1.6pt{\rm R}}}

\newc{\itQQ}[2]{\dspl{\int_{#1}#2\,dz}}
\newc{\mmitQQ}[2]{\dspl{\rule[1mm]{4mm}{.3mm}\hspace{-4.3mm}\int_{#1}~#2~dz}}
\newc{\MmitQQ}[2]{\dspl{\rule[1mm]{4mm}{.3mm}\hspace{-4.3mm}\int_{#1}~#2~d\mu}}

\newc{\MUmitQQ}[3]{\dspl{\rule[1mm]{4mm}{.3mm}\hspace{-4.3mm}\int_{#1}~#2~d#3}}
\newc{\MUitQQ}[3]{\dspl{\int_{#1}~#2~d#3}}

\newc{\mccP}{\mathbb{P}}
\newc{\mccK}{\mathbb{K}}

\newc{\DKTmU}{\mccK(U)}
\newc{\DKTmUold}{(K_U(U)^{-1})^T}

\newc{\myPi}{\mathbf{W}}
\newc{\myIbar}{\bar{\ccI}_1}
\newc{\myIhat}{\hat{\ccI}_1}
\newc{\myIbreve}{\breve{\ccI}_0}

\newc{\mmk}{\mathbf{k}}

\newcommand{\ma}{\mathbf{a}}
\newcommand{\mg}{\mathbf{g}}

\newc{\mfu}{\mathbf{f_u}}
\newc{\mh}{\mathbf{h}}
\newc{\mb}{\mathbf{b}}

\newcommand{\barrl}[2]{\barr{ll}\lefteqn{#1}\hspace{#2}&\\}

\newc{\mN}{\mathbf{N}}
\newc{\mI}{\mathbf{I}}
\newc{\mH}{\mathbf{H}}

\newc{\mk}{\mathbf{k}}
\newc{\mr}{\mathbf{r}}

\newc{\DIAGM}[2]{\left[\barr{ccc}#1&0\ldots&0\\
	\vdots&\ddots&\vdots\\0&\ldots0&#2\earr \right]}
\newc{\DiagM}[2]{\mbox{diag}\left[#1
	\cdots #2 \right]}
\newc{\vVEC}[2]{\left[\barr{c}#1\\
	\vdots\\#2\earr \right]}
\newc{\hVEC}[2]{\left[#1
	\cdots #2 \right]}

\newc{\mq}{\mathbf{q}}
\newcommand{\refexam}[1]{Example~\ref{#1}}

\begin{document}

\vspace*{-.8in}
\begin{center} {\LARGE\em On the Positivity of Weak Solutions to a Class of Cross Diffusion Systems.}

 \end{center}

\vspace{.1in}

\begin{center}

{\sc Dung Le}{\footnote {Department of Mathematics, University of
Texas at San
Antonio, One UTSA Circle, San Antonio, TX 78249. {\tt Email: Dung.Le@utsa.edu}
}}

\end{center}

\begin{abstract}
We establish the positivity of weak (and very weak) solutions to a class of cross diffusion systems which is  inspired by models in mathematical biology/ecology, in particular the Shigesada-Kawasaki-Teramoto (SKT) model in population biology. Examples and counterexamples will be provided to show that our contions are near optimal.  \end{abstract}

\vspace{.2in}

\section{Introduction}\label{intro}
Let $T_0>0$. We consider the following system
\beqno{mabsys0}\left\{\barr{l}u_t=\Div(\ma(u)Du)+\mb(u)Du+ \mg(u)u=0\mbox{ in $Q_0=\Og\times(0,T_0)$},\\u=0 \mbox{ on $\partial\Og \times (0,T_0)$},\\ u(x,0) = u_0(x) \mbox{ on $\Og $}.\earr \right.\eeq

Here, $\ma,\mb,\mg$ are matrices  of size $m\times m$ depending on $u=[u_1,\ldots,u_m]^T\in\RR^m$, The entries of $\ma$ are scalars 
and the entries $b_{ij}$ of the matrix $\mb$ are row vectors (so that  $b_{ij}^T\in \RR^N$). The initial data $u_0$ is a vector value function on $\Og$, a bounded domain with smooth boundary $\partial\Og$ in $\RR^n$.

We are interested in the positivity of (weak and strong) solutions of \mref{mabsys0}. That is, whether the initial data $\psi$ is positive then so are the corresponding (weak) solutions. This property plays an important role in mathematical models  modeling biology/ecology phenomena by evolution systems. Such property was well known for scalar equations and there is a tremendous literature on this topic. However very little is  known for systems like \mref{mabsys0} due to the lack of comparison principles because the classical techniques for scalar equations are not extended easily to systems due to the presence of derivatives in the diffusion parts.

Our first goal is the positivity of (weak) solutions of the {\em quasi-linear} system \mref{mabsys0}. Using a dual argument, we will show that the positivity property of \mref{mabsys0} can be reduced to that of dual {\em linear} systems with smooth coefficients.

We  are then led to the study of positivity property of linear systems with smooth coefficients. This can be done by a change of variables to diagonalize the diffusion parts and then inspection of reaction parts to obtain positive answers.

The results are accompanied with examples showing different ways of transformations in applications. Counterexamples are also provided to show that the results are near optimal in the sense that in some case if our conditions are violated then one can  have either negative or positive answers in some cases.

The notations and main results will be described in \refsec{mainres}. The proof of the duality is presented in \refsec{possec}. We turn to strong solutions of linear systems in \refsec{postrong}. Finally, examples and counterexamples are gathered in \refsec{example}.

\section{Main results}\label{mainres}\eqnoset
We say that $u$ is a very weak solution of \mref{mabsys0} on $Q_0=\Og\times(0,T_0)$ if there are numbers $p_*, q_*\ge 1$ such that for all $T\in(0,T_0)$ we have $u(T)\in L^1(\Og)$ (here, $u(T)(x)=u(x,T)$, and we often drop the letter $x$ in notations if it is clear from the context), $u\in L^{q_*}(Q_0)$ and $Du\in L^{p_*}(Q_0)$ and $u$ satisfies for all  $\fg\in C^1(Q)$, $Q=\Og\times(0,T)$
\beqno{Pwdefdiv}\barrl{\iidx{\Og}{\myprod{u(T),\fg(T)}-\myprod{u_0,\fg(0)}}=}{1cm}&\itQ{\Og \times(0,T)}{[\myprod{u, \fg_t}-\myprod{Du,\ma^T(u)D\fg} -\myprod{u,\Div(\mb^T(u)\fg)}+ \myprod{u,\mg^T(u)\fg}]}.\earr\eeq
We note that the entries $b_{ij}$ of the matrix $\mb$ are row vectors.

Given such an $u$, the main idea of the duality argument is to look for $\fg$ solving
\beqno{keydual}\itQ{\Og \times(0,T)}{[\myprod{u, \fg_t}-\myprod{Du,\ma^T(u)D\fg} -\myprod{u,\Div(\mb^T(u)\fg)}+ \myprod{u,\mg^T(u)\fg}]}=0.\eeq
For such $\fg$ $$\iidx{\Og}{\myprod{u(T),\fg(T)}}=\iidx{\Og}{\myprod{u_0,\fg(0)}}.$$

Through out this work, for any vector function $u$ we will say $u\ge0$ if the components of $u$ are not negative.
Given that $u_0\ge0$, for any $\fg(0)\ge0$ if we can find $\fg$ such that \mref{keydual} holds and $\fg(T)\ge0$ then the above implies that $u(T)\ge0$. We will show that this scheme can be done if some continuity assumptions on the data of the system \mref{mabsys0} are assumed and certain positivity principle holds for the corresponding linear system \mref{keydual}. 

To be precise, suppose that we can find a number $\sg_N\ge1$ such that for all solution $\fg$ of \mref{keydual} and $Q=\Og\times(0,T)$ there is some constant $C$ such that \beqno{dualbound}\|\fg\|_{L^{\sg_N}(Q)}<C.\eeq  Recalling the definition of $p_*,q_*$ in that of weak solutions, we assume that there is a positive number $\sg_N$ such that the following continuity conditions hold, 
\bdes\item[c1)] The map $\ma:L^{q_*}(Q)\to L^{q_1}(Q)$ is continuous for $q_1=\frac{2p_*}{p_*-2}$ if $p_*>2$ and $q_1=\infty$ if $p_*=2$ ;
\item[c2)] the map $\mb:L^{q_*}(Q)\to L^{q_2}(Q)$ is continuous for $q_2=\frac{\sg_N'p_*}{p_*-\sg_N'}$ if $p_*>\sg_N'$ and $q_2=\infty$ if $p_*=\sg_N'$;
\item[c3)] the map $\mg:L^{q_*}(Q)\to L^{q_3}(Q)$ is continuous for $q_3=\frac{\sg_N'q_*}{q_*-\sg_N'}$ if $q_*>\sg_N'$ and $q_3=\infty$ if $q_*=\sg_N'$.
\edes

Our first result is the following dual positivity principle.

\btheo{dualposthm} Let $u$ is a weak solution of \mref{mabsys0}. Suppose that there is a sequence $\{u_n\}$ of sufficiently smooth functions  converges to $u$ in $L^{q_*}(Q_0)$.
For any given $\psi\in C^{1}(\Og)$,  consider the following linear parabolic system with smooth coefficients
\beqno{Psinsysdiv}\left\{\barr{l}W_t=\Div(\ma(u_n)^T DW)-\Div(\mb^T(u_n) W) + \mg(u_n)^T W\mbox{ in $Q_0$},\\ W=0 \mbox{ on $\partial\Og \times (0,T_0)$},\\ W(x,0) = \psi(x) \mbox{ on $\Og $}.\earr \right.\eeq

Assume that  the norms $\|DW_n\|_{L^2(Q)},\; \|W_n\|_{L^{\sg_N}(Q)}$ are uniformly bounded (see \mref{dualbound}) and the continuity conditions c1)--c3) hold. Suppose also that positivity principle holds for \mref{Psinsysdiv}. That is,  strong solutions $W_n$ of \mref{Psinsysdiv} is non-negative for any non-negative initial data $\psi\in C^1(\Og)$.

If $u_0\ge0$ in $\Og$ 
then $u$  is also non-negative in $Q_0$.\etheo

Of course. we can take $\{u_n\}$ to be a sequence of mollifications of $u$. That is, we consider $C^\infty$ functions $\eta(t)$ and $\rg (x)$ whose supports are $(-1,1)$ and $B_1(0)$ and $\|\eta\|_{L^1(\RR)} =\|\rg\|_{L^1(\RR^N)}=1$. Denote $\eta_n(t) = n\eta(t/n)$ and $\rg_n(x)=n^N\rg (x/n)$.
We define $$ u_{n}(t,y)=(\eta_n\fg_n)* u(t,y)=\int_{\RR}\iidx{\RR^N}{\eta_n (s-t)\fg_n(x-y)u(t,x)}\,ds.$$
Then the existence of strong solutions $W_n$ and their uniform estimates assumed in \reftheo{dualposthm} can be established as we know that $u_n\to u$ in $L^r(Q)$ if $u\in L^r(Q)$.

\reftheo{dualposthm} leads us to the investigation of the positivity principle of {\em linear} parabolic systems. We want to investigate the positivity of the following system (compare with \mref{Psinsysdiv})
\beqno{Wsysdiv0}\left\{\barr{l}W_t=\Div(\ma DW)-\Div(\mb W) + \mg W\mbox{ in $Q$},\\ W=0 \mbox{ on $\partial\Og \times (0,T)$},\\ W(x,0) = \psi_0(x) \mbox{ on $\Og $}.\earr \right.\eeq 

Because $u_n$'s are smooth, we will always assume in this section that $\ma,\mb,\mg$ are smooth. 

We will show that 

\btheo{mAmGthm} Let $W$ is a strong solution of \mref{Wsysdiv0} and $\ccJ$ be a differentiable matrix in $(x,t)$ such that \bdes\item[i)] $\ccJ\ma \ccJ^{-1}, \ccJ\mb \ccJ^{-1}$ and $D\ccJ\ccJ^{-1}$ are diagonal; \item[ii)] the  off diagonal entries of $\hat{\mG}=\ccJ\mg\ccJ^{-1}+\ccJ_t\ccJ^{-1}$  are nonnegative.\edes

If $\ccJ(x,0) \fg(x)\ge0$ then $\ccJ(x,t) W\ge0$.
\etheo

In some cases, we have that $\ccJ^{-1}$ is a positive matrix (i.e. its entries are positive). Then $\ccJ W\ge0$ implies also that $W\ge0$. For examples, taking $\ccJ=I$, we have
\bcoro{FTcoro} If $\ma$ is diagonal then the positivity holds if the off-diagonal entries of $\mg$ are nonnegative. \ecoro

A bit more generally, we have

\bcoro{dualposcoro} Consider the system \mref{mabsys0}.
Assume that there is a {\em constant} matrix $\ccJ$ such that $\ccJ \ma(z)\ccJ^{-1}$, $\ccJ \mb(z)\ccJ^{-1}$ are diagonal and the off diagonal entries of $\ccJ \mg(z)\ccJ^{-1}$ are nonnegative for all $z\in \RR^m$.
If $u$ is a solution of \mref{mabsys0} with $(\ccJ^T)^{-1}u_0\ge0$ then $(\ccJ^T)^{-1}u\ge0$.
\ecoro

\reftheo{mAmGthm} and \refcoro{dualposcoro} also  say that if there is a {\em constant} matrix $\ccJ$ such that $\ccJ \ma(z)\ccJ^{-1}$, $\ccJ \mb(z)\ccJ^{-1}$ are diagonal and the off diagonal entries of $\ccJ \mg(z)\ccJ^{-1}$ are nonnegative for all $z\in \RR^m$ then the positive cone of $\RR^m$ is {\em invariant} under  $\ccJ S_1, (\ccJ^T)^{-1}S_2$. Here, $S_1,S_2$ are respectively the solution flows of \mref{Wsysdiv0} and \mref{mabsys0}.

Examples of \reftheo{mAmGthm} showing that its conditions can be relaxed as well as its counterexamples will be presented in \refsec{example}.

\section{Proof of dual positivity principle} \label{possec}\eqnoset

For each $n$ and any given $\psi\in C^{1}(\Og)$ and $T\in(0,T_0)$,  there are strong solutions $\Psi_{n}$ to the following linear parabolic system with smooth coefficients ($\hat{u}_{n}(x,t)=u_{n}(x,T-t)$. Compare with \mref{Psinsysdiv})
\beqno{Psinsysdiv0}\left\{\barr{l}\hat{\Psi}_t=\Div(\ma(\hat{u}_n)^T D\hat{\Psi})-\Div(\mb^T(\hat{u}_n)\hat{\Psi}) + \mg(\hat{u}_n)^T\hat{\Psi}\mbox{ in $Q=\Og\times (0,T)$},\\ \hat{\Psi}=0 \mbox{ on $\partial\Og \times (0,T_0)$},\\ \hat{\Psi}(x,0) = \psi(x) \mbox{ on $\Og $}.\earr \right.\eeq

We first need the following technical lemma. For the sake of convenience, we repeat the continuity conditions c1)--c3) here.

\blemm{duallemma} Suppose that there is a constant $C_0$ such that
\beqno{D2Wn}  \|D\hat{\Psi}_n\|_{L^2(Q)},\; \|\hat{\Psi}_n\|_{L^{\sg_N}(Q)}\le C_0 \quad \forall n,\eeq
and that there is $\sg_N>0$ such that
\bdes\item[i)] For $q_1=\frac{2p_*}{p_*-2}$ if $p_*>2$ and $q_1=\infty$ if $p_*=2$ the map $\ma:L^{q_*}(Q)\to L^{q_1}(Q)$ is continuous;
\item[ii)] For $q_2=\frac{\sg_N'p_*}{p_*-\sg_N'}$ if $p_*>\sg_N'$ and $q_2=\infty$ if $p_*=\sg_N'$ the map $\mb:L^{q_*}(Q)\to L^{q_2}(Q)$ is continuous;
\item[iii)] For $q_3=\frac{\sg_N'q_*}{q_*-\sg_N'}$ if $q_*>\sg_N'$ and $q_3=\infty$ if $q_*=\sg_N'$ the map $\mg:L^{q_*}(Q)\to L^{q_3}(Q)$ is continuous.
\edes

Then \beqno{dualuPsi} \iidx{\Og}{\myprod{u(T),\psi}}=\lim_{n\to\infty}\iidx{\Og}{\myprod{u_0,\hat{\Psi}_n(T)}}.\eeq
\elemm

\bproof By a change of variables $t\to T-t$, the functions $\Psi(x,t)=\hat{\Psi_n}(x,T-t)$ satisfies
\beqno{Psinsys}\left\{\barr{l}\Psi_t+\Div(\ma(u_n)^TD\Psi) -\Div(\mb^T(u_n)\Psi)+ \mg(u_n)^T\Psi=0\mbox{ in $Q:=\Og \times (0,T)$},\\\Psi=0 \mbox{ on $\partial\Og \times (0,T)$},\\\Psi(x,T) = \psi(x) \mbox{ on $\Og $}.\earr \right.\eeq
 
Multiplying the system of $\Psi_n$ by $u$ and integrating over $Q$, we have 
$$\itQ{Q}{[\myprod{u,(\Psi_n)_t}-\myprod{\ma(u_n)Du, D\Psi_n}-\myprod{u,\Div(\mb^T(u_n)\Psi_n)}+\myprod{\mg(u_n)u, \Psi_n}]}=0.$$

From the above equation  and \mref{Pwdefdiv} with $\fg =\Psi_n$ (this is eligible because $\Psi_n$ is a strong solution), using \mref{Psinsys} and the fact that $\Psi_n(T,x)=\psi(x)$, we derive
\beqno{Pwdefz1}\barrl{\iidx{\Og}{[\myprod{u(x.T),\psi}-\myprod{u_0,\Psi_n(x,0)}]}=\itQ{\Og \times(0,T)}{\myprod{[\ma(u_n)-\ma(u)]Du,D\Psi_n}} +}{5cm}
& \itQ{\Og \times(0,T)}{\myprod{[\mb(u_n)-\mb(u)]Du,\Psi_n}}+
\\& \itQ{\Og \times(0,T)}{\myprod{[\mg(u)-\mg(u_n)]u,\Psi_n}}.\earr\eeq

Letting $n\to\infty$, we will see that the integrals on the right hand side tend to zero. Indeed, from the assumption \mref{D2Wn} \beqno{D2Wn0}  \|D\Psi_n\|_{L^2(Q)},\; \|\Psi_n\|_{L^{\sg_N}(Q)}\le C_0 \quad \forall n.\eeq

We consider the first integral in \mref{Pwdefz1}. By \mref{D2Wn0} we have $\|D\Psi_n\|_{L^2(Q)}$'s are bounded uniformly in $n$ so that  we need only to show that $\ma(u_n)Du\to \ma(u)Du$ strongly to 0 in $L^2(Q)$. Let $q=q_1$. that is, $q=\frac{2p_*}{p_*-2}$ if $p_*>2$ and $q=\infty$ if $p_*=2$. By H\"older's inequality,
$$\|[\ma(u_{n})-\ma(u)]Du\|_{L^2(Q)}\le \|\ma(u_{n})-\ma(u)\|_{L^{q}(Q)}\|Du\|_{L^{p_*}(Q)}.$$

As we are assuming in i) that the map $u\to \ma(u)$ is continuous from $L^{q_*}(Q)$ to $L^q(Q)$ and because $u_{n}\to u$ in
$L^{q_*}(Q)$, it is clear that $\ma(u_{n})$ converges to $\ma(u)$ in $L^q(Q)$. Thus, $[\ma(u_{n})-\ma(u)]Du$ converges strongly to 0 in $L^2(Q)$. Thus, the first integral on the right hand side of \mref{Pwdefz1} tends to 0 as $n\to \infty$. 

Similar arguments apply to the second and third integrals on the right hand side of \mref{Pwdefz1} to obtain the same conclusion. Using \mref{D2Wn0}, we see that $\|\Psi_n\|_{L^{\sg_N}(Q)}$'s are uniformly bounded. Thus, if $[\mb(u_n)-\mb(u)]Du$ and $[\mg(u_{n})-\mg(u)]u$ converge strongly to 0 in $L^{{\sg_N}'}(Q)$. With $q_2,q_3$ being defined in ii), iii), by H\"older's inequality, we have
$$\|[\mb(u_{n})-\mb(u)]Du\|_{L^{\sg_N'}(Q)}\le \|\mb(u_{n})-\mb(u)\|_{L^{q_2}(Q)}\|Du\|_{L^{p_*}(Q)},$$
$$\|[\mg(u_{n})-\mg(u)]u\|_{L^{\sg_N'}(Q)}\le \|\mg(u_{n})-\mg(u)\|_{L^{q_3}(Q)}\|u\|_{L^{q_*}(Q)}.$$
These terms goes to zero by the continuity assumptions of $\mb,\mg$. 

We just prove that the right hand side of \mref{Pwdefz1} tends to 0. Thus,
$$\lim_{n\to\infty}\iidx{\Og}{[\myprod{u(x,T),\psi}-\myprod{u_0,\Psi_n(x,0)}]}=0.$$ Since $\Psi_n(x,0)=\hat{\Psi}_n(x,T)$  we have 
$$\iidx{\Og}{\myprod{u(x,T),\psi}}=\lim_{n\to\infty}\iidx{\Og}{\myprod{u_0,\Psi_n(x,0)}}=\lim_{n\to\infty}\iidx{\Og}{\myprod{u_0,\hat{\Psi}_n(x,T)}}.$$

We proved the Lemma. \eproof

We now can present 

{\bf Proof of \reftheo{dualposthm}:}  For any nonnegative $\psi\in C^{1}(\Og)$, let $\hat{\Psi}_n$ be strong solutions of \mref{Psinsysdiv}. By the assumption that the norms $\|DW_n\|_{L^2(Q)},\; \|W_n\|_{L^{\sg_N}(Q)}$ are uniformly bounded \mref{D2Wn} holds and \reflemm{duallemma} applies. By the positity assumption $W_n=\hat{\Psi}_n\ge0$ we have  $$\iidx{\Og\times\{T\}}{u(x,T)\psi(x)}=\lim_{n\to\infty}\iidx{\Og\times\{T\}}{\hat{\Psi}_n(x,T)u_0(x)} \ge0.$$ So, we conclude that $u(x,T)\ge 0$ for all $T\in(0,T_0)$. Hence, $u\ge0$ on $Q$. \eproof

Let us discuss the assumption \mref{D2Wn} and the number $\sg_N$ in the continuity conditions i)--iii) of \reflemm{duallemma}. To this end, we recall the following parabolic Sobolev imbedding inequality, which can be proved easily by using the H\"older and Sobolev imbedding inequalities.

\blemm{parasobolev} Let $r^*=p/N$ if $N>p$ and $r^*$ be any number in $(0,1)$ if $N\le p$. For any sufficiently nonegative smooth functions $g,G$ and any time interval $I$ there is a constant $C$ such that\beqno{paraSobo}\itQ{\Og\times I}{g^{r^*}G^p}\le C\sup_I\left(\iidx{\Og\times\{t\}}{g}\right)^{r^*}\itQ{\Og\times I}{(|DG|^p+G^p)}\eeq   If $G=0$ on the parabolic boundary $\partial\Og\times I$ then the integral of $G^p$ over $\Og\times I$ on the right hand side can be dropped.

Furthermore, if $r<r^*$ then for any $\eg>0$ we can find a constant $C(\eg)$ such that
\beqno{paraSobo1}\itQ{\Og\times I}{g^{r}G^p}\le C\sup_I\left(\iidx{\Og\times\{t\}}{g}\right)^{r}\itQ{\Og\times I}{(\eg|DG|^p+C(\eg)G^p)}\eeq
\elemm

We then have

\blemm{DWlem} Assume that there is a constant $C_0$ such that
\beqno{gintcond0} \sup_{t\in(0,T)}\||\mb|^2+|\mg|\|_{L^{q_0}(\Og\times\{t\})}\le C_0 , \quad \mbox{ for some $q_0>N/2$}.\eeq

Then \mref{D2Wn}  of \reflemm{duallemma} holds for $\sg_N=4/N+2$.

\elemm

\bproof Multiplying the system for $\hat{\Psi}$ by $\hat{\Psi}$ we get for all $t\in(0,T)$ and $Q_t=\Og\times(0,t)$ $$ \iidx{\Og\times\{t\}}{|\hat{\Psi}|^2}+
\itQ{Q_t}{|D\hat{\Psi}|^2}\le \iidx{\Og}{|\psi|^2}+\itQ{Q_t}{(|\mb|\|\hat{\Psi}||D\hat{\Psi}|+|\mg||\hat{\Psi}|^2)}.$$ Applying Young's inequality, we get \beqno{keyH0} 
\iidx{\Og\times\{t\}}{|\hat{\Psi}|^2}+
\itQ{Q_t}{|D\hat{\Psi}|^2}\le \iidx{\Og}{|\psi|^2}+C\itQ{Q_t}{(|\mb^2|+|\mg|)|\hat{\Psi}|^2}.\eeq

By H\"older inequality we estimate the right hand side  by $$\int_0^t \left(\iidx{\Og\times\{\tau\}}{(|\mb|^2+|\mg|)^{q_0}}\right)^\frac{1}{q_0}\left(\iidx{\Og\times\{\tau\}}{|\hat{\Psi}|^{2q_0'}}\right)^\frac{1}{q_0'}d\tau.$$
Since $q_0>N/2$, $2q_0'<2N/(N-2)$. The Sobolev inequality shows that for any $\eg>0$ there is $C(\eg)>0$ such that
$$\left(\iidx{\Og\times\{\tau\}}{|\hat{\Psi}|^{2q_0'}}\right)^\frac{1}{q_0'}\le \eg\iidx{\Og\times\{\tau \}}{|D\hat{\Psi}|^2}+C(\eg)\iidx{\Og\times\{t\}}{|\hat{\Psi}|^2}.$$
This and \mref{gintcond0} imply that $$\itQ{\Og\times(0,t)}{(|\mb|^2+|g|)|\hat{\Psi}|^2}\le C(C_0)\eg\itQ{Q_t}{|D\hat{\Psi}|^2}+C(\eg,C_0)\itQ{Q_t}{|\hat{\Psi}|^2}.$$

Using this in \mref{keyH0} with $\eg$ sufficiently small, we get
$$\iidx{\Og\times\{t\}}{|\hat{\Psi}|^2}+
\itQ{Q_t}{|D\hat{\Psi}|^2}\le \iidx{\Og}{|\psi|^2}+C(C_0)\itQ{Q_t}{|\hat{\Psi}|^2}.$$ From this Gr\"onwall inequality for $\|\hat{\Psi}\|_{L^2(\Og\times\{t\})}^2$ we get $\sup_t\|\hat{\Psi}\|_{L^2(\Og\times\{t\})}^2\le C(C_0)$. We also obtain from \mref{keyH0} and the above estimate that $$\itQ{Q_T}{|D\hat{\Psi}|^2}\le C(C_0,\|\psi\|_{L^{2}(\Og)}).$$

Applying \reflemm{parasobolev} with $g=|\hat{\Psi}|$, $G=\hat{\Psi}$ and $p=2$, we conclude that
$$\|\hat{\Psi}\|_{L^{\sg_N}(Q_T})\le C(C_0,\|\psi\|_{L^{2}(\Og)}),\quad 
\sg_N=4/N+2.$$

The lemma is proved. \eproof

\subsection{The special case $\ma=P_u$} We now consider the case when there is a map $P:\RR^m\to\RR^m$ such that $\ma =P_u$. We consider the following system
\beqno{Psys}\left\{\barr{l}u_t=\Delta(P(u))+\mb(u)Du+ \mg(u)u=0\mbox{ in $Q_0$},\\u=0 \mbox{ on $\partial\Og \times (0,T_0)$},\\ u(x,0) = \psi(x) \mbox{ on $\Og $}.\earr \right.\eeq

We say that a $u$ is a very weak solution of \mref{Psys} if  for all $T\in(0,T_0)$ $u(T)\in L^1(\Og)$, $u\in L^{q_*}(Q_0)$ and $Du\in L^{p_*}(Q_0)$ and $u$ satisfies for all admissible $\fg$ with $\fg_t\in C^1(Q)$ and $\Delta \fg \in C^2(Q)$
\beqno{PwdefdivP}\barrl{\iidx{\Og}{\myprod{u(T),\fg(T)}-\myprod{u_0,\fg(0)}}=}{1cm}&\itQ{\Og \times(0,T)}{[\myprod{u, \fg_t}+\myprod{P(u),\Delta\fg} -\myprod{u,\Div(\mb^T(u)\fg)}+ \myprod{u,\mg^T(u)\fg}]}.\earr\eeq

We can write $P(u)=\bar{\ma}(u)u$ with \beqno{madef} \bar{\ma}(u)=\int_0^1 P_u(su)\,ds.\eeq The system \mref{Psys} is then written in divergence form \mref{mabsys0}. By \mref{PwdefdivP}
The dual system \mref{Psinsysdiv} now can be
the linear parabolic system
\beqno{PWsysdiv}\left\{\barr{l}\hat{\Psi}_t=\bar{\ma}^T(\hat{u}_n)\Delta\hat{\Psi}-\Div(\mb^T(\hat{u}_n) \hat{\Psi}) + \mg^T(\hat{u}_n)\hat{\Psi}\mbox{ in $Q_0$},\\ \hat{\Psi}=0 \mbox{ on $\partial\Og \times (0,T_0)$},\\ \hat{\Psi}(x,0) = \psi(x) \mbox{ on $\Og $}\earr \right.\eeq 
Testing the system with $\Delta\hat{\Psi}_n$ and arguing as in \reflemm{DWlem}, we can prove that (see also \cite[Lemma 3.5 and Remark 3.4]{lduni}) if
\beqno{gintcond} \sup_{t\in(0,T)}\||\mg|^2\|_{L^{q_0}(\Og\times\{t\})}\le C_0 , \quad \mbox{for some $q_0>N/2$}\eeq then there is a constant $C(\|D\psi\|_{L^2(\Og)})$ such that for all $n$ \beqno{D2Wnz}  \sup_{(0,T)}\|D\Psi_n\|_{L^2(\Og)},\; \|\Delta\Psi_n\|_{L^2(Q)}\le C(\|D\psi\|_{L^2(\Og)}).\eeq

We can use \reflemm{parasobolev} with $g=|D\hat{\Psi}_n|^2$, $|G=D\hat{\Psi}_n$ and $p=2$ to see that $D\hat{\Psi}_n\in L^q(Q)$ where $q\in(1,\infty)$ if $N=2$ and $q=\frac4N+2$. As $D\hat{\Psi}_n\in L^2(\Og)$, $\hat{\Psi}_n\in L^\frac{2N}{N-2}(\Og)$. We apply \reflemm{parasobolev} again with $g=|D\hat{\Psi}_n|^\frac{2N}{N-2}$ and $G=\hat{\Psi}_n$ and $p=q=\frac{4}{N}+2$ to see that
\beqno{Psisg}\|\hat{\Psi}_n\|_{L^{\sg_N}(Q)}\le C,\eeq
where \beqno{newsg}\sg_N=\left\{\barr{ll}\mbox{any number in $(1,\infty)$} & \mbox{if $N=2$},\\ \mbox{any number in $(1,6+\frac{10}{3})$} &\mbox{if $N=3$},\\\frac{2(4+2N)}{N-2}+\frac{4}{N}+2 &\mbox{if $N>3$}.
\earr\right.\eeq

\blemm{duallemma1} Assume \mref{gintcond0}. Let $\sg_N$ be defined by \mref{newsg}. Suppose that the continuity conditions hold with c1) (or i) of \reflemm{duallemma}) being replaced by \bdes\item[c1')] For $q_1=\frac{2q_*}{q_*-2}$  if  the map $\bar{\ma}:L^{q_*}(Q)\to L^{q_1}(Q)$ is continuous.
\edes

Then the conclusion \mref{dualuPsi} of \reflemm{duallemma} still holds with $\hat{\Psi}_n$ being strong solutions of \mref{PWsysdiv}.
\elemm

The lemma is proved by the same argument. Again, we define $\Psi_n(x,t)=\hat{\Psi}_n(x, T-t)$.
As before, we multiply the equation \mref{PWsysdiv} with $u$ and use \mref{PwdefdivP} and $P(u)=\bar{\ma}(u)u$ to see that \mref{Pwdefz1} now is \beqno{PPwdefz1}\barrl{\iidx{\Og}{[\myprod{u(T),\psi}-\myprod{u_0,\Psi_n(0)}]}=\itQ{\Og \times(0,T)}{\myprod{[\bar{\ma}(u_n)-\bar{\ma}(u)]u,\Delta\Psi_n}} +}{5cm}
& \itQ{\Og \times(0,T)}{\myprod{[\mb(u_n)-\mb(u)]Du,\Psi_n}}+
\\& \itQ{\Og \times(0,T)}{\myprod{[\mg(u)-\mg(u_n)]u,\Psi_n}}.\earr\eeq

Using the uniform bound of $\|\Delta\Psi_n\|_{L^2(Q)}$ in \mref{D2Wnz}, if $u\in L^{q_*}(Q)$ then by H\"older's inequality we easily see that $|u||\Delta \Psi_n| \in L^p(Q)$ with $p=\frac{2q_*}{q_*+2}$. As $u_n\to u$ in $L^{q_*}(Q)$ and $p'=q_1=\frac{2q_*}{q_*-2}$, under our assumption c1') the first integral on the right hand side of \mref{PPwdefz1} tends to 0. 

Again, we have the same conlusion of \reflemm{duallemma} with the assumptions of \reflemm{duallemma1}.

\section{Positivity of the dual system:} \label{postrong}\eqnoset
The discussion in the previous section leads us to the investigation of positivity principle of linear parabolic systems \mref{Psinsysdiv0}. In this section, we consider the following linear  system (denoting $W=\hat{\Psi}_n$ and $\ma,\mg$ by $\ma^T, \mg^T$ for simplicity)
\beqno{Wsysdiv}\left\{\barr{l}W_t=\Div(\ma DW)-\Div(\mb W) + \mg W\mbox{ in $Q$},\\ W=0 \mbox{ on $\partial\Og \times (0,T)$},\\ W(x,0) = \psi_0(x) \mbox{ on $\Og $}.\earr \right.\eeq 

We will always assume in this section that $\ma,\mb,\mg$ are smooth as $u_n$ is. In the rest of this paper, we will denote by $W$ a strong solution of the above system. 
Let $\ccJ=\ccJ(x,t)$ be a differentiable and invertible $m\times m$ matrix and $v=\ccJ W$. We multiply the system of $W$ with $\ccJ$ and show that $v$ satisfies a system of the form
\beqno{vsysdiv}\left\{\barr{l}v_t=\Div(\mA Dv)-\Div(\mB_1 v) +\mB_2 Dv+ \mG v\mbox{ in $Q$},\\ v=0 \mbox{ on $\partial\Og \times (0,T)$},\\ v(x,0) = \ccJ(x,0)\psi_0(x) \mbox{ on $\Og $}.\earr \right.\eeq
where $\mA, \mB$ and $\mG$ are matrices in terms of $\ma,\mb,\mg$ and $\ccJ$. Of course, the entries of $\mA, \mG$ are scalar functions and those of $\mB_1, \mB_2$ are vectors in $\RR^N$.

We first have this positivity result for \mref{vsysdiv}.
\bprop{mAmG} Let $v$ be a weak solution of \mref{vsysdiv}. Suppose that $\mA, \mB_1, \mB_2$ are diagonal matrices and that there are constants $\llg_1>0, C_0$ such that $\myprod{\mA \zeta,\zeta}\ge\llg_1|\zeta|^2$ and
\beqno{mGcond}\myprod{-\mG v,v^-}\le C_0|v^-|^2.\eeq If $\ccJ(x,0)\psi_0(x)\ge0$ for all $x\in\Og$ then $v\ge0$.
\eprop

\bproof   Because $Dv_i^+Dv_i^-=v_i^+Dv_i^+=0$, we know that $$\myprod{Dv,Dv^-}=\myprod{Dv^+-Dv^-, Dv^-}=-\myprod{Dv^-, Dv^-}$$ and (the vectors)  $\myprod{v,Dv^-}=\myprod{Dv,v^-}=-\myprod{ v^-,Dv^-}$. We test the system of \mref{vsysdiv} of $v$ by $v^-$. That is, we multiply the $i$-th equation with $v_i^-$ and sum the results. Using integration by parts and the assumption that $\mA, \mB_1, \mB_2$ are diagonal matrices, we have for $\mB=\mB_1+\mB_2$ (whose entries are vectors)
$$-\frac{d}{dt}\iidx{\Og}{|v^-|^2}-\iidx{\Og}{\myprod{\mA Dv^-,Dv^-}}= -\iidx{\Og}{\myprod{\mB v^-,Dv^-}}+ \iidx{\Og}{\myprod{\mG v,v^-}}.$$

Because $\myprod{\mA \zeta,\zeta}\ge\llg_1|\zeta|^2$, for $\zeta=Dv^-$ we have $\myprod{\mA Dv^-,Dv^-}\ge \llg_1|Dv^-|^2$. Using this for the integral of $\myprod{\mA Dv^-,Dv^-}$ and applying Young's inequality to the intgral of $\myprod{\mB v^-,Dv^-}$, we get $$\frac{d}{dt}\iidx{\Og}{|v^-|^2}+\llg_1\iidx{\Og}{|Dv^-|^2}\le C_1|v^-|^2- \iidx{\Og}{\myprod{\mG v,v^-}}.$$

Using \mref{mGcond}, we then have
\beqno{gronwall}\frac{d}{dt}\iidx{\Og}{|v^-|^2}+\llg_1\iidx{\Og}{|Dv^-|^2}\le  C_2\iidx{\Og}{|v^-|^2}.\eeq

As $v(0)=\ccJ(0)\psi_0\ge0$, we then have $v^-(0)=0$. The above inequality yields a Gr\"onwall inequality for $\|v^-\|_{L^2(\Og\times\{t\})}^2$ and implies that $v^-=0$. Thus, $v\ge0$.
\eproof 

\brem{posDir} The result still holds if $v$ satisfies non-homogeneous Dirichlet boundary condition $v=\fg$ on $\partial\Og\times(0,T)$ with $\fg\ge0$. In this case, $v^-=0$ on $\partial\Og\times(0,T)$ because $\fg\ge0$ and hence there is no boundary integrals in the integration by parts. The calculation and argument remain the same.

Similarly, the result also holds for homogeneous Neumann boundary condition.
\erem

\newc{\mf}{\mathbf{f}}
\brem{frem} We can introduce a term $\mf$ to the right hand side of the system \mref{vsysdiv}. The result and proof will be the same if we assume that $\ccJ\mf\ge0$.
\erem

\blemm{mGlem} The condition \mref{mGcond} holds if the off-diagonal entries of $\mG$ are nonnegative. \elemm
\bproof
We have $\myprod{\mG v,v^-}=\sum_{i,j}g_{ij}v_iv^-_j$. If $i=j$ then $v_iv^-_j=-(v^-_i)^2$. Otherwise, $v_iv^-_j=v^+_iv^-_j-v^-_iv^-_j\ge -v^-_iv^-_j$. 
Thus, if $g_{i,j}\ge0$ for $i\ne j$ then $\myprod{\mG v,v^-}\ge -C\sum_{i,j}v_i^-v^-_j$. By H\"older's inequality $\myprod{-\mG v,v^-}\le C\sum_{i,j}v_i^-v^-_j\le C|v^-|^2$. \eproof

From \refprop{mAmG} and \reflemm{mGlem}, we will prove that (this is \reftheo{mAmGthm})
\bprop{mAmGx} Assume that there is an invertible matrix $\ccJ(x,t)$ such that \bdes\item[i)] $\ccJ\ma \ccJ^{-1}, \ccJ\mb \ccJ^{-1}$ and $D\ccJ\ccJ^{-1}$ are diagonal; \item[ii)] the  off diagonal entries of $\hat{\mG}=\ccJ\mg\ccJ^{-1}+\ccJ_t\ccJ^{-1}$  are nonnegative.\edes

If $\ccJ(x,0) \psi_0(x)\ge0$ then $\ccJ(x,t) W\ge0$.
\eprop

The assertion is also true for the dual system \mref{PWsysdiv} of \mref{Psys}
\beqno{aPusys} W_t=\ma\Delta W -\Div(\mb W)+\mg W.\eeq

\brem{ccJdiag} Concerning the diagonality of $D\ccJ\ccJ^{-1}$, we also
observe that $D\ccJ^{-1}=-\ccJ^{-1}D\ccJ \ccJ^{-1}$ so that $D\ccJ\ccJ^{-1}=-\ccJ D(\ccJ^{-1})$.

\erem

We will prove \refcoro{mAmGx} by considering the cases $\ccJ$ depends only on $t$ first and then $\ccJ$ depends also on $x$ and showing that $v=\ccJ W$ satisfies the conditions of \refprop{mAmG}.

Before this, we can present

{\bf Proof of \refcoro{dualposcoro}:} Let  $u$ be a solution of \mref{mabsys0} and $W_n$'s be solutions of the corresponding \mref{Psinsysdiv}. If we take $\ccJ$ to be a constant matrix then by \reflemm{duallemma}
$$\iidx{\Og}{\myprod{u(T),\psi_0}} = \lim \iidx{\Og}{\myprod{u_0,W_n(T)}}.$$
This also implies
\beqno{x1}\iidx{\Og}{\myprod{(\ccJ^T)^{-1}u(T),\ccJ\psi_0}} = \lim \iidx{\Og}{\myprod{(\ccJ^T)^{-1}u_0,\ccJ W_n(T)}}.\eeq

Because $\ccJ \ma(z)\ccJ^{-1}, \ccJ \mb(z)\ccJ^{-1}$ are diagonal and the off diagonal entries of $\ccJ \mg(z)\ccJ^{-1}$ are nonnegative for all $z\in \RR^m$, the same assertions are true for   $\ccJ \ma^T(u_n)\ccJ^{-1}, \ccJ \mb^T(u_n)\ccJ^{-1}$ and  $\ccJ \mg^T(u_n)\ccJ^{-1}$. We then can apply \reftheo{mAmGthm} to \mref{Psinsysdiv}  to see that if $\ccJ\psi_0\ge0$ then $\ccJ W_n(T)\ge0$.  As $(\ccJ^T)^{-1}u_0\ge0$, \mref{x1} implies $$\iidx{\Og}{\myprod{(\ccJ^T)^{-1}u(T),\ccJ\psi_0}}\ge0$$  For any $C^1$ vector function $\fg\ge0$ we take $\psi_0=\ccJ^{-1}\fg$ (then $\ccJ\psi_0\ge0$) and see that
$$\iidx{\Og}{\myprod{(\ccJ^T)^{-1}u(T),\fg}}\ge0.$$ This implies $(\ccJ^T)^{-1}u(T)\ge0$. This is true for all $T\in(0,T_0)$. We proved the assertion. \eproof

\subsection{$\ccJ(x,t)=\ccJ(t)$ is independent of $x$}
Let $\ccJ$ be an invertible  matrix independent of $x$ and $v=\ccJ W$. Multiply the system \mref{Wsysdiv} of $W$ with $\ccJ$,
because $v_t=\ccJ_t W+\ccJ W_t=\ccJ W_t+\ccJ_t\ccJ^{-1}v$. Using the system for $W$, the matrices $\mA$, $\mB_i$ and $\mG$ in \mref{vsysdiv} are \beqno{mABG}\mA=\ccJ\ma \ccJ^{-1},\;\mB_1=\ccJ\mb \ccJ^{-1},\; \mB_2=0 \mbox{ and }  \mG = \ccJ\mg\ccJ^{-1}+\ccJ_t\ccJ^{-1}.\eeq

Thus, if $\ccJ \ma\ccJ^{-1},\ccJ\mb\ccJ^{-1}$ are diagonal then \refprop{mAmG} and and \reflemm{mGlem} apply to give \refprop{mAmGx}.

\bexam{crosssys} Consider a cross diffusion system of 2 equations $u_t=\Div(\ma Du)+\mg u$ where $\ma=[\ag_{ij}]$ and $\mg=[g_{ij}]$ are constant $2\times2$ matrices
$$\ma=\left[\barr{cc}\ag_{11}&\ag_{12}\\\ag_{21}&\ag_{22}\earr\right], \; \mg=\left[\barr{cc}g_{11}&g_{12}\\g_{21}&g_{22}\earr\right].$$
Assume that $\ma$ has  positive eigenvalues $d_1\ne d_2$  and two linear independent eigenvectors $e_1, e_2$ with positive components. This implies that there are constants $a,b,c,d>0$ such that $\ccJ\ma\ccJ^{-1}=\mbox{diag}[d_1,d_2]$ with $$\ccJ^{-1}=\left[\barr{cc}a&b\\c&d\earr\right],\; \ccJ=\frac{1}{ad-bc}\left[\barr{cc}d&-b\\-c&a\earr\right].$$
The off diagonal entries of $\ccJ\mg\ccJ^{-1}$ are
$$\frac{a^2g_{12}-b^2g_{21} -ab(g_{11}-g_{22})}{ad-bc},\;\frac{-c^2g_{12}+d^2g_{21} +cd(g_{11}-g_{22})}{ad-bc}.$$

If these numbers are nonnegative then \refcoro{mAmGx} (of course, $D\ccJ \ccJ^{-1}=\ccJ_t\ccJ^{-1}=0$ here) implies that if $\ccJ\psi_0\ge0$ then $v=\ccJ W\ge0$. Since $\ccJ^{-1}$ is positive, we have $W\ge0$.
\eexam

\subsection{$\ccJ$ depends on $x,t$} \eqnoset We consider the case $\ccJ$ depends also on $x$ and establish \refprop{mAmG} again.
As before we let $v=\ccJ W$. We multiply the system of $W$ by $\ccJ$ and compute $v_t, Dv$ to obtain a similar system for $v$.
We still have $v_t=\ccJ W_t+\ccJ_t\ccJ^{-1}v$  and need to compute $\ccJ \Div(\ma DW)$.

First of all, $Dv=\ccJ DW+D\ccJ W$ so that $\ccJ DW= Dv -D\ccJ W$. 
We also have $\Div(\ccJ\ma DW)=\ccJ\Div(\ma DW)+D\ccJ\ma DW$, Thus, for $\mA=\ccJ \ma \ccJ^{-1}$ $$\barr{lll}\ccJ\Div(\ma DW)&=&\Div(\ccJ\ma DW)-D\ccJ\ma DW=\Div(\mA \ccJ DW)-D\ccJ\ma\ccJ^{-1} \ccJ DW\\
&=&\Div(\mA (Dv-D\ccJ W))-D\ccJ\ma\ccJ^{-1} (Dv-D\ccJ W)\\&=&
\Div(\mA Dv)-\Div(\mA D\ccJ W)-D\ccJ\ma \ccJ^{-1}Dv+D\ccJ\ma\ccJ^{-1} D\ccJ W\\&=&
\Div(\mA Dv)-\Div(\mA D\ccJ\ccJ^{-1} v)-D\ccJ\ma \ccJ^{-1}Dv+D\ccJ\ma\ccJ^{-1} D\ccJ\ccJ^{-1}v
\earr$$

Similarly,
$$\ccJ\Div(\mb W)=\Div(\ccJ\mb W)-D\ccJ\mb W=\Div(\ccJ \mb \ccJ^{-1} v)-D\ccJ\mb \ccJ^{-1}v.$$

We see that $v$ satisfies \mref{vsysdiv} $v_t=\Div(\mA Dv)-\Div(\mB_1 v) +\mB_2 Dv+ \mG v$, where
$$\mA=\ccJ\ma\ccJ^{-1},\; \mB_1=\mA D\ccJ\ccJ^{-1}+\ccJ \mb \ccJ^{-1},\; \mB_2=-D\ccJ\ma \ccJ^{-1},$$
$$\mG=\ccJ\mg^T\ccJ^{-1}+D\ccJ\mb \ccJ^{-1}+\ccJ_t\ccJ^{-1}+D\ccJ\ma\ccJ^{-1} D\ccJ\ccJ^{-1}.$$

We see that $\mB_2=-D\ccJ\ccJ^{-1}\mA$. Therefore, the matrices $\mA, \mB_1$ and $\mB_2$ are diagonal if $\ccJ\ma \ccJ^{-1}, \ccJ\mb \ccJ^{-1}$ and $D\ccJ\ccJ^{-1}$ are diagonal.  

Furthermore, $D\ccJ\ma\ccJ^{-1} D\ccJ\ccJ^{-1}=D\ccJ\ccJ^{-1}\mA D\ccJ\ccJ^{-1}$ and $D\ccJ\mb\ccJ^{-1}=D\ccJ\ccJ^{-1}\ccJ\mb\ccJ^{-1}$ are diagonal. Thus, the off diagonal entries of $\mG$ are those of
$$\hat{\mG}=\ccJ\mg^T\ccJ^{-1}+\ccJ_t\ccJ^{-1}.$$

Thus, $v=\ccJ W$ solves a system satisfying the conditions of \refprop{mAmG}.

\subsection{The special case $\ma=P_u$}

We consider the dual system \mref{aPusys} of \mref{Psys} $$W_t=\ma\Delta W -\Div(\mb W)+\mg W.$$
We have $\Div(\ma DW)=\ma \Delta W+ D\ma DW$ so that $$\barr{lll}\ccJ\ma\Delta W&=&\ccJ\Div(\ma DW)-\ccJ D\ma DW=\ccJ\Div(\ma DW)-\ccJ D\ma \ccJ^{-1}\ccJ DW\\
&=& \ccJ\Div(\ma DW)-\ccJ D\ma \ccJ^{-1}(Dv -D\ccJ W)\\
&=& \ccJ\Div(\ma DW)-\ccJ D\ma \ccJ^{-1}Dv +\ccJ D\ma \ccJ^{-1}D\ccJ \ccJ^{-1} v.\earr$$

The term $\ccJ\Div(\ma DW)$ will be treated as before. We find that $v$ satisfies  $$v_t=\Div(\mA Dv)-\Div(\mB_1 v) +\mB_2 Dv+ \mG v,$$ where
$$\mA=\ccJ\ma\ccJ^{-1},\; \mB_1=\mA D\ccJ\ccJ^{-1}+\ccJ \mb \ccJ^{-1},\; \mB_2=-D\ccJ\ma \ccJ^{-1}-\ccJ D\ma \ccJ^{-1},$$ $$\bar{\mG}=\ccJ\mg^T\ccJ^{-1}+D\ccJ\mb \ccJ^{-1}+\ccJ D\ma \ccJ^{-1}D\ccJ \ccJ^{-1}+D\ccJ\ma\ccJ^{-1} D\ccJ\ccJ^{-1}.$$

We now look at the new matrix factors $\ccJ D\ma \ccJ^{-1}$ of $Dv$ and $\ccJ D\ma \ccJ^{-1}D\ccJ \ccJ^{-1}$ of  $v$. We have $\ccJ\ma =\mA\ccJ$ so that $D(\ccJ\ma) =D(\mA\ccJ)$.
This implies $$\ccJ D\ma =D\mA\ccJ+\mA D\ccJ -D\ccJ \ma \mbox{ or }
\ccJ D\ma \ccJ^{-1}=D\mA+\mA D\ccJ\ccJ^{-1} -D\ccJ \ma\ccJ^{-1}.$$

Because $D\mA$, $\mA D\ccJ\ccJ^{-1}$ and $D\ccJ \ma\ccJ^{-1}=D\ccJ \ccJ^{-1}\mA$ are diagonal we conclude that $\ccJ D\ma \ccJ^{-1}$ is diagonal and so is $\mB_2$. Since $D\ccJ \ccJ^{-1}$ is diagonal, $\ccJ D\ma \ccJ^{-1}D\ccJ \ccJ^{-1}$  is also diagonal.
We then see that the previous argument remains.

\section{Examples}\label{example}\eqnoset

In this section, we present some ways of choosing $\ccJ$ to apply \reftheo{mAmGthm} to  the following system (compare with \mref{Psinsysdiv})
\beqno{Wsysdiv0u}\left\{\barr{l}u_t=\Div(\ma Du)-\Div(\mb u) + \mg u\mbox{ in $Q$},\\ u=0 \mbox{ on $\partial\Og \times (0,T)$},\\ u(x,0) = u_0(x) \mbox{ on $\Og $}.\earr \right.\eeq 

The diagonality  $\ccJ\ma\ccJ^{-1},\ccJ\mb\ccJ^{-1}, D\ccJ \ccJ^{-1}$ and the nonnegativity of the off diagonal entries of $\ccJ\mg\ccJ^{-1}$ are the key conditions and we will also provide counterexamples when its conditions are violated.

The first choice of $\ccJ$ is $\ccJ=\mbox{diag}[e^{\cg_i}]\ccC$ where $\cg_i$'s are functions in $x,t$ and $\ccC$ is an invertible matrix in $t$. 

First of all, it is clear that if $\ccC\ma\ccC^{-1}$  and $\ccC\mb\ccC^{-1}$  are diagonal then so are $\ccJ\ma\ccJ^{-1}$ 
 and $\ccJ\mb\ccJ^{-1}$. Secondly, since $D\ccJ=\mbox{diag}[ e^{\cg_i}D\cg_i]\ccC$ and $\ccJ^{-1}=\ccC^{-1}\mbox{diag}[e^{-\cg_i}]$, $D\ccJ\ccJ^{-1}=\mbox{diag}[D\cg_i]$ is diagonal.

Finally, we look at the signs of the off diagonal entries of $\mG=\ccJ\mg\ccJ^{-1}+\ccJ_t\ccJ^{-1}$.  Because $\ccJ_t=\mbox{diag}[(\cg_i)_t e^{\cg_i}]\ccC+\mbox{diag}[e^{\cg_i}]\ccC_t$ so that $$\ccJ \mg^T\ccJ^{-1}=\mbox{diag}[e^{\cg_i}]\ccC\mg^T\ccC^{-1}\mbox{diag}[e^{-\cg_i}],\;\ccJ_t\ccJ^{-1}=\mbox{diag}[(\cg_i)_t]+\mbox{diag}[e^{\cg_i}]\ccC_t\ccC^{-1}\mbox{diag}[e^{-\cg_i}].$$

The off-diagonal entries of $\mG$ then are those of $\mbox{diag}[e^{\cg_i}](\ccC \mg^T\ccC^{-1}+\ccC_t\ccC^{-1})\mbox{diag}[e^{-\cg_i}]$ and their signs are those of $\ccC \mg^T\ccC^{-1}+\ccC_t\ccC^{-1}$.

Assume that $\ma=I$. Let  $k$ be a positive constant and $K$ be a constant matrix. Assume that $K$ is commuting with $\mg^T$ and the off diagonal entries of $K$ are positive. Let $\ccC=e^{kKt}$. Then $\hat{\mG}=\mg^T+kK$. If  we choose $k$ large then the off diagonal entries of $\hat{\mG}$ are positive.

If $\ma,\mg$ are commute then they can be diagonalized simultaneously by some matrix $\ccC$. We then assume that $\ccC_t\ccC^{-1}$ has nonnegative off diagonal entries.

Assume that $\ma=\ccC^{-1}\ag\ccC$ for some diagonal matrix $\ag=\mbox{diag}[\ag_i]$ where $\ag_i$'s are functions in $x,t$. We see that
$\mA=\ccJ\ma\ccJ^{-1}=\mbox{diag}[e^{\cg_i}]\ag\mbox{diag}[e^{-\cg_i}]$ is diagonal.

\reftheo{mAmGthm} asserts that  $\ccJ W\ge0$ if $\ccJ(x,0)\psi_0\ge0$. Of course, if $\ccJ^{-1}$ is a nonnegative matrix (whose entries are nonnegative) then $W=\ccJ^{-1}\ccJ W\ge0$. In the sequel, we will be interested in this case. For simplicity, we will consider the case when $\ma$ is diagonal and $\mb=0$ and $\ccJ$ depends only on $t$. The diagonality of $\ccJ\ma\ccJ^{-1}, \ccJ\mb\ccJ^{-1},D\ccJ\ccJ^{-1}$ is obvious. We will concentrate on the nonnegativity of the off diagonal entries of $\mG$.

Let $n_i$, $i=1,\ldots, k$, be integers such that $\sum n_i=m$. We consider the case $\ma=\DiagM{d_1(w) \mI_1}{d_k(w) \mI_k}$ where $d_i$'s are positive functions on $\RR^m$ and $\mI_i$'s are identity matrices of size $n_i\times n_i$. If $\ccJ_i$'s are $n_i\times n_i$ matrices and $\ccJ:=\DiagM{\ccJ_1}{\ccJ_k}$, then $\ma$,  $\ccJ$ commute (blocks by blocks)
and it is clear that $\ccJ\ma\ccJ^{-1}=\ma$ is diagonal.

We state the following obvious result for later references.
\bcoro{posKG} Assume that there are $n_i\times n_i$ matrices $\ccJ_i$'s such that for $\ccJ:=\DiagM{\ccJ_1}{\ccJ_k}$
\bdes\item[i)] $\DiagM{\ccJ_1^{-1}}{\ccJ_k^{-1}}$ is positive;
\item[ii)] the off diagonal entries of $\ccJ\mg\ccJ^{-1}+\DiagM{(\ccJ_1)_t\ccJ_1^{-1}}{(\ccJ_k)_t\ccJ_k^{-1}}$ are nonnegative.
\edes

If $\ccJ(0)\psi_0\ge0$ then $W\ge0$.
\ecoro

In particular, we take $\ccJ=\DiagM{e^{K_1t}}{e^{K_kt}}$ where $K_i$'s are constant matrices of size $n_i\times n_i$ in $t$. We
then have that $\ccJ_t=\DiagM{K_1e^{K_1 t}}{K_ke^{K_k t}}$ and $\ccJ^{-1}=\DiagM{e^{-K_1 t}}{e^{-K_k t}}$. Since $\ccJ_t\ccJ^{-1}=\DiagM{K_1}{K_k}$ and hence $$\hat{\mG}=\ccJ\mg\ccJ^{-1}+\DiagM{K_1}{K_k}.$$

We will see that for an appropriate choice of $K_i$'s, the off-diagonal entries of $\hat{\mG}$ are nonnegative even those of $\mg$ are not.
We then have the following result.

\bcoro{posKG1} Assume that
\bdes\item[i)] the matrices $e^{-K_i t}$'s have positive entries;
\item[ii)] the off diagonal entries of $\ccJ\mg\ccJ^{-1}+\DiagM{K_1}{K_k}$ are nonnegative.
\edes

If $\psi_0\ge0$ then $W\ge0$.
\ecoro

\bproof We take $\ccJ=\DiagM{e^{K_1t}}{e^{K_kt}}$. Since $\hat{\mG}$ has positive entries by ii), we have that $v=\ccJ W\ge0$ (as $\ccJ(0)=\mI$ and $\psi_0\ge0$ we have $\ccJ(0)\psi_0\ge0$). Because $\ccJ^{-1}(t)=\DiagM{e^{-K_1 t}}{e^{-K_i t}}$ is a positive matrix, by i), we also have $W=\ccJ^{-1}v\ge0$.
\eproof

In the following examples we will show that the assumption that diagonal entries of $\mg$ are positive (the condition ii) of \refcoro{posKG} and \refcoro{posKG1})  could be relaxed. We consider the system of two equations
\beqno{sysex1} W_t = \Div(\ma DW)+\mg W \mbox{ where  $\ma=d(W)\mI$ for some function $d$}.\eeq

We will apply \refcoro{posKG} to this system. The matrix $\ccJ$ will be of the form $e^{K(t)}$ where $K(t)=\ag(t)\mI+\mN(t)$ where $\mN(t)$ is either {\em nilpotent} or {\em idempotent}. Of course, the examples also apply to the case when $\ma$ can be diagonalised by $\ccJ$.

To begin, we have the following remark for later reference.

\brem{genK} If $K(t)=\ag(t)\mI +\mN(t)$ with $\mN(t)$ is nilpotent (i.e., $\mN^m(t)=0$ for some integer $m$) and $\ccJ=e^{K(t)}=e^{\ag(t)\mI}e^{\mN(t)}$ then $\ccJ^{-1}=e^{-\ag(t)\mI}e^{-\mN(t)}$ and $\ccJ_t=\ag'(t)e^{\ag(t)\mI}e^{\mN(t)}+e^{\ag(t)\mI}(e^{\mN(t)})_t$. Note that $e^{\mN(t)}=\sum_{i=0}^{m-1}\frac{N^i(t)}{i!}$. So that
$$\ccJ_t\ccJ^{-1}=\ag'(t)\mI+(e^{\mN(t)})_te^{-\mN(t)}=\ag'(t)\mI+\left(\sum_{i=0}^{m-1}\frac{\mN^i(t)}{i!}\right)_t\left(\sum_{i=0}^{m-1}\frac{(-\mN(t))^i}{i!}\right)$$
In particular, if $m=2$ then $\ccJ_t\ccJ^{-1}=\ag'(t)\mI+\mN'(t)(\mI-\mN(t))$.

Similarly, if $K(t)=\ag(t)\mI +\mN(t)$ with $\mN(t)$ is idempotent (i.e., $\mN^m(t)=\mN(t)$ for every integer $m\ge2$)  then we note that $e^{\mN(t)}=\sum_{i=0}^{m-1}\frac{N^i(t)}{i!}=\mI+(e-1)\mN(t)$ and $e^{-\mN(t)}=\mI+(e^{-1}-1)\mN(t)$. So that for $\ccJ=e^{K(t)}$
$$\ccJ_t\ccJ^{-1}=\ag'(t)\mI+(e^{\mN(t)})'e^{-\mN(t)}=\ag'(t)\mI+(e-1)N'(t)(\mI+(e^{-1}-1)\mN(t)).$$

We note that a $2\times2$ matrix $\mN$ is idempotent iff $\mN= \left[\barr{cccc}a&b\\c&1-a\earr\right]$ with $a^2+bc=a$.
\erem

\bexam{ex1} Let $a,b,\eg$ be functions with $\eg(t)\ge0$ and $\int_0^\infty\eg(s)dx<\infty$ we take 
$$\mg=\left[\barr{cc}a(t)&-\eg(t)\\0&a(t)\earr\right],\;  K(t)=\left[\barr{cc}t&b(t)\\0&t\earr\right].$$ We observe that $G, K$ commute so that $\hat{\mG}=\mg+\ccJ_t\ccJ^{-1}$. Because $$\ccJ=e^{K(t)}=e^t\left[\barr{cc}1&b(t)\\0&1\earr\right],\;\ccJ_t=e^t\left[\barr{cc}1&b(t)\\0&1\earr\right]+e^t\left[\barr{cc}0&b'(t)\\0&0\earr\right].$$ Thus, $\ccJ_t\ccJ^{-1}=\mI+\left[\barr{cc}0&b'(t)\\0&0\earr\right]$ and $\mg+\ccJ_t\ccJ^{-1}=\left[\barr{cc}a(t)+1&-\eg(t)+b'(t)\\0&a(t)+1\earr\right]$.
Hence, the off diagonal entries of $\hat{\mG}$ are nonnegative if $b'(t)\ge\eg(t)$ and the entries of $\ccJ^{-1}$ are nonnegative if $b(t)\le0$. This is the case if we choose $b(t)=\int_0^t\eg(s)ds-\int_0^\infty\eg(s)ds$. \eexam

\bexam{balone0} If $\mg=[g_{ij}]$ is a $2\times2$ matrix with $g_{ij}$ being functions in $t$.   Let $\ccJ=e^\mN$ with $\mN=\left[\barr{cc}0&b(t)\\0&0\earr\right]$. Then $\ccJ^{-1}=\mI-\mN$.  Denoting $\mq(b)=g_{21}b^2+(g_{11}-g_{22})b-g_{12}$, a straightforward calculation shows that
$$\ccJ\mg\ccJ^{-1}=\left[\barr{cc}g_{11}+bg_{21}&-\mq(b)\\g_{21}&-bg_{21}+g_{22}\earr\right],\; \ccJ_t\ccJ^{-1}=\left[\barr{cc}0&b'(t)\\0&0\earr\right].$$
Thus, the off diagonal entries of $\ccJ\mg\ccJ^{-1}+\ccJ_t\ccJ^{-1}$ will be nonnegative if   $g_{21}\ge0$ and $b'=\mq(b)$.

Assume that $g_{21}=0$. The solution of $b'=\mq(b)$ is
\beqno{sol1}b(t)=e^{\int_0^t(g_{11}(\tau)-g_{22}(\tau))ds}\left(b(0)-\int_0^te^{-\int_0^s(g_{11}(\tau)-g_{22}(\tau))d\tau}g_{12}(s)ds\right).\eeq

Thus, if $g_{12}<0$ and $$\int_0^\infty e^{-\int_0^s(g_{11}(\tau)-g_{22}(\tau))d\tau}g_{12}(s)ds=b(0)$$
is a (negative) finite number then $b(t)<0$ for all $t\in(0,\infty)$. So, there is still a nonnegative solution $W$ if the initial data $\psi_0$ of $W$ satisfies $\psi_{0,1}(x,0)+b(0)\psi_{0,2}(x,0), \psi_{0,2}(x,0)\ge0$. Indeed, this condition means $\ccJ(0)\psi_0(0)\ge0$. As the off diagonal entries of $\hat{\mG}$ are nonnegative, $\ccJ(t)W(t)\ge0$.  Because $b(t)<0$  so that $\ccJ^{-1}(t)$ is positive and therefore $W(t)\ge0$ for $t>0$.
\eexam

Let us consider the case when $g_{ij}$'s are constants with $g_{21}>0$ and $g_{12}<0$.
Since $g_{21}>0$, the dynamics of the solution of $b'=\mq(b)$ is well known. If $\delta:=(g_{11}-g_{22})^2+4g_{21}g_{12}<0$ then for any given $b(0)$, $b(t)$ is increasing and blows up to $\infty$ in finite time. If $\delta>0$ then the quadratics $\mq(b)$ has two distinct roots $b_1<b_2<0$ or $0<b_1<b_2$. The first case happens if $g_{11}>g_{22}$. In this case, if $b(0)\le b_2$ then $b(t)$ is either increasing or decreasing to $b_1$ and   $b(t)\le0$ for all $t\in (0,\infty)$.

Accordingly, we can conclude that 
\bcoro{nondiag}
If $g_{21}>0$ and $g_{12}<0$ then there is a smooth function $b$ on $[0,t_0)$, $t_0$ is finite if $\delta=(g_{11}-g_{22})^2+4g_{21}g_{12}<0$ and $t_0=\infty$ if $\delta>0$,  such that if  $\psi_{0,1}(x,0)+b(0)\psi_{0,2}(x,0), \psi_{0,2}(x,0)\ge0$  then $W_{1}(x,t)+b(t)W_{2}(x,t), W_{2}(x,t)\ge0$ for all $t\in(0,t_0)$. 

In particular, assume that $\psi_{0,1}(x,0)+b(0)\psi_{0,2}(x,0), \psi_{0,2}(x,0)\ge0$ for some constant $b(0)<0$. If $\delta<0$ then there is a finite $t_0$ such that $W(t)\ge0$ for all $t\in (0,t_0)$. If $\delta>0$, $g_{11}>g_{22}$ and $b(0)\le b_2$, the larger root of $\mq(b)$, then $W(t)\ge0$ for all $t\in(0,\infty)$.
\ecoro

The Corollary is a consequence of \refcoro{posKG} with $\ccJ=e^\mN=\mI+\mN$. Indeed, the first assertion comes from the fact that $\ccJ(0)\psi_0\ge0$ so that $\ccJ(t)W(t)\ge0$ while the second one comes from $\ccJ^{-1}(t)=\mI-\mN$ is a positive matrix as long as $b(t)\le0$.

\bexam{calone}
As another example, we will take $\mN=\left[\barr{cccc}1&0\\c(t)&0\earr\right]$ and $\ccJ=e^\mN$ then $\mN$ is idempotent. A direct calculation shows that
$$\ccJ=e^{\mN(t)}=\left[\barr{cccc}e&0\\(e-1)c(t)&1\earr\right],\; \ccJ^{-1}=e^{-\mN(t)}=\left[\barr{cccc}e^{-1}&0\\(e^{-1}-1)c(t)&1\earr\right],$$ $$\ccJ_t\ccJ^{-1}=\left[\barr{cccc}0&0\\(1-e^{-1})c'(t)&0\earr\right].$$

If $\mg=[g_{ij}]$, a $2\times2$ matrix, then $\ccJ\mg\ccJ^{-1}$ is $$\left[ \begin {array}{cc}  g_{{11}}+(1-e)g_{{12}}c &eg_{{12}}\\ \noalign{\medskip}
-(1-e^{-1})\mr & \left( e-1 \right) cg_{{12}}+g_{{22}}
\end {array} \right] 
$$ with $\mr=-(e-1)g_{12}c^2-(g_{11}-g_{22})c-\frac{1}{e-1}g_{21}$.
Therefore, the off diagonal entries of $\ccJ\mg\ccJ^{-1}+\ccJ_t\ccJ$ are nonnegative if $$g_{12}\ge0\mbox{ and }-(1-e^{-1})\mr 
+(1-e^{-1})c'=0.$$
This is to say $g_{12}\ge0$ and
$c'=\mr$.

Consider the case $g_{12}=0$. We then need a function $c$ such that
$$c'=-(g_{11}-g_{22})c-\frac{1}{e-1}g_{21}.$$ The solution of this is
$$c(t)=e^{\int_0^t(g_{11}-g_{22})ds}\left(c(0)-\frac{1}{e-1}\int_0^te^{-\int_0^s(g_{11}-g_{22})d\tau}g_{21}(s)ds\right).$$

Thus, if $g_{21}\le0$ and $(e-1)c(0)-\int_0^\infty e^{-\int_0^s(g_{11}-g_{22})d\tau}g_{21}(s)ds\le0$ then $c(t)\le0$ for all $t$ and $\ccJ^{-1}=e^{-\mN(t)}$ is positive. If the initial data $\psi_0$ of $W$ satisfies  $\ccJ(0)\psi_0\ge0$ then, as the off diagonal entries of $\hat{\mG}$ are nonnegative, $\ccJ(t)W(t)\ge0$.  Because  $\ccJ^{-1}(t)$ is positive (as $c(t)\le0$) so that $W(t)\ge0$.
\eexam

\bexam{bthenc}
We now consider the case $g_{ij}$ are constants and both $g_{12}, g_{21}<0$. We first use $\ccJ=e^{\mN}$ with $\mN=\left[\barr{cc}0&b(t)\\0&0\earr\right]$, a nilpotent matrix. Let $b$ solves $b'=\mq$ with $\mq$ being introduced in \refexam{balone0}. We have $$b'=g_{21}b^2+(g_{11}-g_{22})b-g_{12}.$$

Since $g_{12},g_{21}<0$, if $g_{11}-g_{22}<0$ the quadratics $\mq$ has two negative roots $b_1<b_2<0$. This equation has a solution $b\le0$ on $(0,\infty)$ if its initial value $b(0)\le b_2$. We also note that either $b(t)$ blows up to $-\infty$ or $$\lim_{t\to\infty}b(t)=b_2=\frac{-(g_{11}-g_{22})-\sqrt{(g_{11}-g_{22})^2+4g_{12}g_{21}}}{2g_{21}}.$$ 

We then have that $v=\ccJ W$ solves $v_t=\Div(\ma Dv)+\hat{\mg} v$ with $$\hat{\mg}=\left[\barr{cc}\hat{g}_{11}&0\\\hat{g}_{21}&\hat{g}_{22}\earr\right], $$
where  $\hat{g}_{21}=g_{21}$,  $\hat{g}_{11}=g_{11}+bg_{21}$ and $\hat{g}_{22}=-bg_{21}+g_{22}$.

Next, we let $\hat{\ccJ}=e^{\hat{\mN}}$ with $\hat{\mN}=\left[\barr{cc}1&0\\c(t)&0\earr\right]$, an idempotent matrix, and (see the case $\hat{g}_{12}=0$ in \refexam{calone}) $$c(t)=e^{\int_0^t(\hat{g}_{11}-\hat{g}_{22})ds}\left(c(0)-\frac{1}{e-1}\int_0^te^{-\int_0^s(\hat{g}_{11}-\hat{g}_{22})d\tau}g_{21}(s)ds\right).$$

As $g_{21}<0$, if $c(0)-\frac{1}{e-1}\int_0^\infty e^{-\int_0^s(\hat{g}_{11}-\hat{g}_{22})d\tau}g_{21}(s)ds\le0$ then $c(t)<0$ for all $t\in(0,\infty)$. Note that $\hat{g}_{11}-\hat{g}_{22}=g_{11}-g_{22}+2b(t)g_{21}\to g_{11}-g_{22}+2b_2g_{21}=-\sqrt{(g_{11}-g_{22})^2+4g_{12}g_{21}}$. Thus, $\int_0^\infty e^{-\int_0^s(\hat{g}_{11}-\hat{g}_{22})d\tau}ds=\infty$ and $c(t)$ becomes positive in finite time.

We see that $w=\hat{\ccJ} v$ solves $w_t=\Div(\ma Dw)+\tilde{\mg}w$ where $\tilde{\mg}=[\tilde{g}_{ij}]$ with $\tilde{g}_{12} =\tilde{g}_{21}=0$. Thus, if $\hat{\ccJ}(0)v(0)\ge0$ then $\hat{\ccJ}(t)v(t)\ge0$ for all $t\in(0,\infty)$. Whenever  $c(t)<0$,  $\hat{\ccJ}^{-1}=e^{-\mN}$ is positive and thus $v\ge0$. This implies  $v=\ccJ W\ge0$. Since $\ccJ^{-1}$ is positive as $b(t)\le0$, we also have $W(t)\ge0$ for $t$ such that both $b(t), c(t)$ exist and are negative. We see that this can be asserted only for $t$ in some finite time interval.
\eexam

\bexam{cthenb} We consider the case $g_{12}, g_{21}<0$ again. We use first $\hat{\mN}=\left[\barr{cccc}1&0\\c(t)&0\earr\right]$ and find $c(t)$ such that $c'=\mr$ where $\mr$ is the quadratics introduced in \refexam{calone} and we then consider
$$c'=-(e-1)g_{12}c^2-(g_{11}-g_{22})c-\frac{1}{e-1}g_{21}.$$

Since $-(e-1)g_{12}>0$ and $g_{12}g_{21}>0$, if $(g_{11}-g_{22})^2-4g_{12}g_{21}>0$ then the quadratics $\mr$ has two roots $c_1,c_2$ and they are of the same sign. If $g_{11}-g_{22}<0$ then $c_1<c_2<0$ and there is $c(t)\le0$ for all $t\ge0$ if $c(0)<c_2$. Moreover, $$\lim_{t\to\infty}c(t)=c_1=\frac{(g_{11}-g_{22})-\sqrt{(g_{11}-g_{22})^2-4g_{12}g_{21}}}{-2(e-1)g_{12}}.$$ 

We then have that $v=\hat{\ccJ}W$ satisfies $v_t=\Div(\ma Dv)+\hat{\mg} v$ with $$\hat{\mg}=\left[\barr{cccc}\hat{g}_{11}&\hat{g}_{12}\\0&\hat{g}_{22}\earr\right],$$
with $\hat{g}_{12}=eg_{12}$, $\hat{g}_{11}=g_{11}+(1-e)g_{12}c$ and $\hat{g}_{22}=(e-1)g_{12}c+g_{22}$.

We then apply $\ccJ=e^{\mN}$ with $\mN=\left[\barr{cccc}0&b(t)\\0&0\earr\right]$ to the system of $v$. Choose $b$ such that $$b(t)=e^{\int_0^t(\hat{g}_{11}(\tau)-\hat{g}_{22}(\tau))ds}\left(b(0)-e\int_0^te^{-\int_0^s(\hat{g}_{11}(\tau)-\hat{g}_{22}(\tau))d\tau}g_{12}(s)ds\right).$$
with $b(0)-e\int_0^\infty e^{-\int_0^s(\hat{g}_{11}(\tau)-\hat{g}_{22}(\tau))d\tau}g_{12}(s)ds\le0$ then $b(t)\le0$ for all $t\ge0$. Note that $\hat{g}_{11}-\hat{g}_{22}=g_{11}-g_{22}+2(1-e)g_{12}c(t)\to g_{11}-g_{22}+2(1-e)g_{12}c_1$. We see that $\hat{g}_{11}-\hat{g}_{22}<0$ and
$$\lim_{t\to\infty}\hat{g}_{11}-\hat{g}_{22}=2(g_{11}-g_{22})-\sqrt{(g_{11}-g_{22})^2-4g_{12}g_{21}}<0.$$
Therefore, $-\int_0^s(\hat{g}_{11}-\hat{g}_{22})d\tau>0$ and $\int_0^\infty e^{-\int_0^s(\hat{g}_{11}(\tau)-\hat{g}_{22}(\tau))d\tau}ds=\infty$. We see that $c(t)$ becomes positive in finite time.

We then have that $w=\ccJ v$ solves $w_t=\Div(\ma Dw)+\tilde{g}w$ for some diagonal matrix $\tilde{g}$. Thus, if $\ccJ(0)v(0)\ge0$ then $\ccJ(t) v(t)\ge0$. Because $b\le0$, $\ccJ^{-1}=\mI-\mN$ is positive so that $v(t)\ge0$. But $W=\hat{\ccJ}^{-1}v$ and $\hat{\ccJ}^{-1}=\mI+(e^{-1}-1)\hat{\mN}$ is a positive matrix (whenever $c(t)\le0$), we can only assert that $W\ge0$ in finite time because $c(t)$ becomes positive in finite time.

\eexam

\brem{earem} In the examples, we can replace $\mN$ by $\ag(t)\mI +\mN$
but this will only add $\ag'(t)$ to the diagonal entries of $\hat{\mG}=\ccJ\mg\ccJ^{-1}+\ccJ_t\ccJ^{-1}$. Thus, thí will make no effects to the analysis and we gain nothing new.
\erem

\brem{otherchoice} Other choices of $\ccJ=e^{\mN}$ with $\mN$ being idempotent are $\mN=\left[\barr{cccc}1&b(t)\\0&0\earr\right]$ and $\mN=\left[\barr{cccc}0&b(t)\\0&1\earr\right]$. However, by inspection, these choices will reduce the matrix $\mg$ to $\hat{\mg}$ with $\hat{g}_{12}=0$ and $\hat{g}_{21}$ is of the same sign of $g_{21}$. The analysis will lead to similar conclusions.
\erem

{\bf Counterexamples:} One would expect that $W\ge0$ globally if $g_{12}, g_{21}<0$. That is $W(x,t)\ge0$ for $(x,t)\in\Og\times(0,\infty)$. The following examples are inconclusive. Let $h=\sin x\sin y$ and $\Og=(0,\pi)\times(0,\pi)$. Define $W_1=1+he^t$ and $W_2=1-he^t$. We have
$\Delta W_1=-2he^t,\; \Delta W_2=2he^t,\; he^t=\frac12(W_1-W_2)$.
This shows that $W=[W_1,W_2]^T$ solves $W_t=\Div(\ma DW)+\mg W$ with
$$\ma=\left[\barr{cccc}1&0\\0&1\earr\right],\; \mg=\frac12\left[\barr{cccc}3&-3\\-3&3\earr\right].$$

We see that the initial data $\psi_0=[1+h,1-h]^T\ge0$ but $W$ can become negative when $t$ is large. Note that $W=[1,1]^T$ on the lateral boundary, a non-homogeneous positive Dirichlet boundary condition (see \refrem{posDir}). We also note that the off diagonal entries of $\mg$ are both negative in this case.

However, define $W_1=1+he^{-t}$ and $W_2=1-he^{-t}$. We have
$(W_1)_t=-he^{-t}$, $(W_2)_t=he^{-t}$, $\Delta W_1=-2he^{-t},\; \Delta W_2=2he^{-t},$ and $he^{-t}=\frac12(W_1-W_2)$.
This shows that $W=[W_1,W_2]^T$ solves $W_t=\Div(\ma DW)+\mg W$ again with
$$\ma=\left[\barr{cccc}1&0\\0&1\earr\right],\; \mg=\frac12\left[\barr{cccc}1&-1\\-1&1\earr\right].$$

The off diagonal entries of $\mg$ are negative. However, for any $t>0$ we see that $W(t)\ge0$.

{\bf Different and complex eigenvalues:} With $u_1=1+he^t$ and $u_2=1-he^t$, this also provides a counterexample when $\ma=\left[\barr{cccc}a_1&0\\0&a_2\earr\right]$ with $a_1\ne a_2$. Indeed, we have $u_t=\Div(\ma Du)+\mg u$ with
$$\ma=\left[\barr{cccc}1&0\\0&\frac12\earr\right],\; \mg=\left[\barr{cccc}\frac32&-\frac32\\-1&1
\earr\right].$$

In the same way, let us consider the case when $\ma$ is not real diagonalizable. Again, let $h=\sin x\sin y$ and $\Og=(0,\pi)\times(0,\pi)$. Define $v_1=1+he^t$ and $v_2=1-he^t$. For any $a,b$
we see that $a\Delta v_1-b\Delta v_2=2(b-a)he^t$ and $b\Delta v_1+a\Delta v_2=2(a-b)he^t$. Since $he^t=\frac12(v_1-v_2)$, $v=[v_1,v_2]^T$ satisfies a non-homogeneous positive Diriclet boundary condition $v=[1,1]^T$ and solves $$v_t=\Div(\ma Dv)+\mG v \mbox{ with }
\ma=\frac12\left[\barr{cccc}a&-b\\b&a\earr\right],\; \mG=\frac12\left[\barr{cccc}1+a-b&b-a-1\\b-a-1&a-b+1\earr\right].$$

If $a>0$ then $\ma$ is elliptic. We choose  $b$ sufficiently large to see that the off diagonal entries of $\mG$ are positive but $v$ changes its sign when $t$ is large.

\newc{\pmat}[4]
{\left[
	\begin{array}{c | c}
		#1 & #2 \\
		\hline
		#3 & #4
	\end{array}
	\right]
}

Of course, the above examples can be extended to systems of three or more equations. For example, in block form

$$\ma = \pmat{\ma_1}{0}{0^T}{\ma_2},\; \mg = \pmat{G_1}{X_1}{X_2^T}{G_2},$$ 
where $\ma_1,G_1$ are $2\times 2$  matrices, $\ma_2,G_2$ are $1\times 1$  matrices and $X_1,X_2$ are vectors of size $2\times 1$.

Accordingly, we take $ \ccJ=\pmat{\ccJ_1}{0}{0^T}{\ccJ_2}$ with $\ccJ_2=[1]$ and $\ccJ_1=e^{\mN}$ with $\mN$ being nilpotent or idempotent as in the previous examples.

Clearly, $\ccJ\ma \ccJ^{-1}$ is diagonal and
$$\ccJ\pmat{G_1}{X_1}{X_2^T}{G_2}\ccJ^{-1}=\pmat{\ccJ_1G_1\ccJ_1^{-1}}{\ccJ_1X_1\ccJ_2^{-1}}{\ccJ_2X_2^T\ccJ_1^{-1}}{\ccJ_2G_2\ccJ_2^{-1}}=\pmat{\ccJ_1G_1\ccJ_1^{-1}}{\ccJ_1X_1}{X_2^T\ccJ_1^{-1}}{G_2}.$$ The earlier arguments show that some off diagonal entries of $\mg$ are allowed to be negative and the positivity property still holds.

 For example, $\mN= \left[\barr{cccc}0&b(t)\\0&0\earr\right]$. Then we can find $b(t)\le0$ so that $\ccJ^{-1}$ is positive. If $\ccJ_1 X_1\ge0$ and $X_2^T\ccJ_1^{-1}\ge0$ then the off diagonal entries of $\mG$ are nonnegative so that we can assert the same conclusions. We leave the details to the readers.

\bibliographystyle{plain}

\end{document}